\documentclass[12pt]{amsart}
\usepackage{latexsym}
\usepackage{amsmath,amsfonts,amssymb,enumerate,multicol,tikz}
\usepackage{epsfig}
\usepackage{pstricks}
\usepackage{hyperref}
\usepackage[utf8]{inputenc}

\usepackage{mathrsfs}

\def\square{\Box}

\begin{document}
\title{The notion of space in Grothendieck: from schemes to a geometry of forms} 
\author{John Alexander Cruz Morales \\
Universidad Nacional de Colombia}

\date{}

\maketitle

\begin{abstract}

In this essay we give a general picture about the evolution of Grohendieck's ideas regarding the notion of space. Starting with his fundamental work in algebraic geometry, where he introduces schemes and toposes as generalizations of classical notions of spaces, passing through tame topology and ending with the formulation of a geometry of forms, we show how the ideas of Grothendieck evolved from pure mathematical considerations to physical and philosophical questions about the nature and structure of space and its mathematical models. 
\end{abstract}

\section{Introduction} Alexander Grothendieck is one of most influential mathematicians of the last century and it is not difficult to argue that is one of the most influential through the history of mathematics. The reach and depth of Grothendieck's work have been extensively discussed in many places, see for instance \cite{car,loc,zal} and the bibliography therein. Our goal is more modest since we do not pretend to cover fully Grothendieck's ideas but just focus on a specific, but central part of his work, namely the notion space\footnote{Since this text is thought for a general audience including mathematicians, philosophers, philosophers of mathematics and, in general, anyone interested in the history and philosophy of the notion of space in the grothendieckean work, we will avoid the technicalities and will privilege a conceptual approach.}. \\ 

It is very clear, even in a superficial reading of Grothendieck's work, that the search of an adequate definition of the notion of space is a leitmotif for Grothendieck, from the very beginning of his mathematical career when he was interested in nuclear spaces to the last mathematical reflections in the 80's of the last century when he tried to develop a geometry of forms as a new foundations for topology and the physical space\footnote{It is known that during his retirement in Lasserre since 1991 until his death, Grothendieck wrote extensively about different topics and it would be plausible to think that one of these topics was the notion of space. However, this is a mere speculation, so for the purpose of this text, we do not take into account that period.}, passing through his monumental work in algebraic geometry. \\

The question {\it What is a space?} has been one of the main questions the human being has asked since the ancient times. Many philosophers, natural philosophers\footnote{We will use the expression {\it natural philosopher} as it was understood in 17th and 18th century.} and mathematicians have tried to give an answer\footnote{There is no point to make an extensive list of the people who have been interested in the problem of the space but this list should include people like Plato, Aristotle, Archimedes, Leibniz, Newton, Riemann, Einstein, just to mention some outstanding examples. Clearly Grothendieck follows this line.} and even though we have now a better understanding of the problem we are still far from a definitive answer. The question is hard since it involves not only mathematics but also physics and philosophy. This was cleverly noted by Grothendieck in the 80's\footnote{Of course, Grothendieck was not the first one who notes this. However, we want to emphasize that Grothendieck realized the importance of philosophy and physics as complements of the mathematical approach to the notion of space, since traditionally the studies on Grothendieck's work have not paid attention to this aspect. On the contrary, it is usually believed that Grothendieck was not interested at all in the physical and philosophical aspects of the notion of space, and this is partially true, if we only focus on Grothendieck's ideas before the 80's. Nevertheless, his mind changed a lot, regarding physical and philosophical issues, and this is one of the points we want to discuss in this essay.} and we will return to this point later in this text. \\

In fact, in \cite{grot1} Grothendieck writes about the notion of space: \\

``{\it La notion d'``espace'' est sans doute une des plus anciennes en mathématique. Elle est
si fondamentale dans notre appréhension ``géométrique'' du monde, qu'elle est restée plus
ou moins tacite pendant plus de deux millénaires. C'est au cours du si\`{e}cle écoulé seulement
que cette notion a fini, progressivement, par se détacher de l'emprise tyrannique de la perception immédiate (d'un seul et m\^{e}me ``espace'' qui nous entoure), et de sa théorisation
traditionnelle (“euclidienne”), pour acquérir son autonomie et sa dynamique propres. De
nos jours, elle fait partie des quelques notions les plus universellement et les plus couramment utilisées en mathématique, famili\`{e}re sans doute \`{a} tout mathématicien sans exception.
Notion protiforme d'ailleurs s'il en fut, aux cents et mille visages, selon le type de structures qu'on incorpore \`{a} ces espaces, depuis les plus riches de toutes (telles les vénérables
structures ``euclidiennes'', ou les structures ``affines'' et ``projectives'' , ou encore les structures ``algébriques'' des ``variétés'' de m\^{e}me nom, qui les généralisent et qui assouplissent)
jusqu'aux plus dépouillées : celles o\`{u} tout élément d'information ``quantitatif'' quel qu'il
soit semble disparu sans retour, et o\`{u} ne subsistent plus que la quintessence qualitative
de la notion de ``proximité'' ou de celle de ``limite'', et la version la plus élusive de
l'intuition de la forme (dit ``topologique'')}''. \\

Following Deligne \cite{del} we can think of a space as something for which {\it localization} makes sense. In this approach the key word is {\it localization}. Thus, the question {\it What is a space?} might be replaced by {\it What does localization mean?}. Along this essay we will try to show how Grothendieck approached this question. The basic idea of localization, the one Deligne had in mind, is exemplified by the following situation\footnote{This example is due to Deligne who wrote it to the author in \cite{del}}: Let $X$ be a topological space. In order to define a continuous function $f: X \rightarrow \mathbb{R}$, it suffices to define it on opens subsets $\mathcal{U}_i$ covering $X$ and to check the agreement on two by two intersections. This way the notion of continuous function is defined by localization on $X$.  \\

In the example above there are two important points that play an important role in the grothendieckean notion of space. The first one is the absence of points, replaced by another thing, in the example by open sets, and the second one is that those open sets cover the whole thing. It has been remarked in other places (see \cite{car}) how problematic the notion of point could be in order to understand what a space is. One of the main features in the several generalizations Grothendieck proposes for the notion of space is that for Grothendieck the space is not made by points. Points are just marks that can be put on the space but not its fundamental elements. On the other hand, the idea of covering is essential in order to understand the relation between the local and the global aspects of the space and it also will be important in the discussion between the continuum and the discrete in Grothendieck's geometry of forms as we will see later.  

\subsection{A brief biographical note}

We will give a very short review of Grothendieck's life. More detailed accounts can be found in \cite{car, schar, zal, jac}\footnote{We include this brief note in order to give the reader an idea of Grothendieck's particular life and shed light on how his personality as mathematician, and more general as human being, was formed.}. Alexander Grothendieck was born in Berlin on March 28 1928. Son of Alexander Schapiro and Hanka Grothendieck both active political militants. Due to the activities of his parents, Grothendieck lived nearby Hamburg\footnote{The years in Hamburg had a big impression in Grothendieck. We were informed by Winfried Scharlau \cite{schar1} that Grothendieck visited Hamburg around 2006. He knew this by some local people that informed him on the visit. In some sense, it seems that Grothendieck was recollecting his own steps.} under the tutelage of the family of a protestant priest until 1939 when he met his family\footnote{Grothendieck's parents took part in the Spanish civil war enrolled in the anarchist militias.}. However, his father was taken prisoner in the concentration camp Le Vernet and then deported to Auschwitz where he died and Alexander and his mother were sent to the camp Rieucros. This compromise of his parents with political issues will be very influential in Grothendieck's point of view, both mathematical and non-mathematical. \\

From 1945 to 1948 Grothendieck studied mathematics at Montpellier. In Récoltes et Semailles \cite{grot1}
he recollects the years in Montpellier and how in solitude developed a measure theory close to the one developed by Lebesgue many years ago. This episode marked his life as mathematician. As Grothendieck wrote (see \cite{grot1}) thanks to this he learnt the importance of the solitude in the work of a mathematician. A solitude that accompanied him his whole life beyond mathematics. \\

In 1949 Grothendieck arrived to Paris to Cartan's seminar at the École normale supérieure. This was an important turn on Grothendieck's mathematical education, since he moved from the border (Montpellier) to the center of French mathematics (Paris). In Paris and Nancy, where he wrote his doctoral dissertation under Jean Dieudonné and Laurent Schwartz, Grothendieck could show his enormous mathematical talent. After this, he spent time as postdoctoral researcher in Sao Paulo and Kansas. \\

During his time in Kansas in 1955 Grothendieck's work took an important turn since he left analysis and started a new mathematical road in (algebraic) geometry. It is in this stay when he wrote the main ideas of his work in abelian  categories\footnote{This work is known as {\it T\^{o}hoku} because of the journal where it was published.} and also he started to work in his version of the Riemann-Roch theorem. Thus, his monumental reconstruction of algebraic geometry started and, even more, his quest of an adequate definition of the notion of space really started in an explicit way. \\

Since 1958 until 1970 Grothendieck occupied a position at IHES. With Dieudonné as his scribe, Grothendieck wrote the giant treatise {\it Elements of Algebraic Geometry} (EGA) and ran his celebrated {\it Seminar of Algebraic Geometry} (SGA). Grothendieck was positioned as the main expert in algebraic geometry in the world and won the Fields medal in 1966 for his fundamental work. During this time he introduced many important concepts (schemes, toposes among others) in order to carry out his program of providing a vast generalization of the notion of space.\\

Despite the success of his seminar and the fruitful years at the IHES, Grothendieck resigned his position in 1970 arguing ethical reasons because the founding the IHES was receiving from the military agency\footnote{This is the universal accepted reason on Grothendieck's resignation to the IHES. However, some authors (see \cite{ruelle}, for instance) believe that there could be deeper reasons behind Grothendieck's demission.}. After being involved in ecological activities in the group {\it Survivre et vivre} he moved back to Montpellier and took one position as professor. It is argued that Grothendieck left mathematics after his departing of IHES since he ceased any mathematical publication but this is far from being true. Certainly Grothendieck abandoned the establishment but this does not mean that he abandoned mathematics. In fact, many important Grothendieck's ideas were developed during the late 70's and the beginning of the 80's when he was in Montpellier. \\

In order to have an idea of Grothendieck's mathematical interests during the 80's, it is good to see \cite{grot2}. This influential text\footnote{It is known the story on how Voevodsky, being a student in Moscow, started to learn French just to read Grothendieck's Esquisse. Certainly this reading was very important in Voevodsky's career.} contains the main ideas of Grothendieck after his reconstruction of the algebraic geometry. There we can find discussions about moduli spaces of Riemann surfaces and the Teichm\"{u}ller tower, tame topology (another approach to a generalization of the notion of space), children drawings, anabelian geometry. This text condensed many of the mathematical thoughts written in manuscripts that widely circulated between the mathematical community. \\

It is interesting to note, however, that some ideas concerning the geometry of forms that will be important for us, in the sense that they constitute an important grothendieckean approach to the notion of space and its relations with physics and philosophy, written in an relatively unknown manuscript \cite{mont}, are not discussed at the Esquisse. These ideas were really important for Grothendieck as we can see from this extract of a letter sent to Tsuji on July 4th 1986: \\

``{\it Excusez-moi d'avoir laissé passer quelques jours avant de répondre \`{a} votre lettre précédente
- depuis pr\`{e}s d'un mois je suis lancé sur un nouveau th\`{e}me mathématique. Des nouveaux
fondements de la topologie, dans un esprit tr\`{e}s différent de celui de la topologie générale
- j'ai envie d'appeler cette nouvelle topologie ``géométrie des formes'' ou ``analysis situs'',
et il est bien possible que ce sera l\`{a} le premier texte mathématique que je publierai, apr\`{e}s
mon départ de 1970 (et avant les divers textes prévus, que j'annonce dans l'introduction \`{a}
Récoltes et Semailles)...}'' \\

Aside from the mathematical texts, Grothendieck wrote in the 80's other general texts of great interest. His long text Récoltes et Semailles is a jewel in the whole sense of the word. There Grothendieck not only explored his past as mathematician but also discussed the intricate aspects of mathematical creativity and the ethics of a mathematician. In La clef des songes, Grothendieck analysed his dreams and arrived to a personal discovering of God. \\

In 1991 Grothendieck decided to move to Lasserre and only kept contact with few people. During this time he wrote extensive texts about different topics, but it is still a mystery the exact content of the thousands of pages he wrote in his isolation. He passed away on November 13 2014 at Saint-Girons, leaving behind him a big legacy for the future generations.

\section{A walk around algebraic geometry}

In this section we will discuss two generalizations of the notion of space that Grothendieck proposed as part of his program of reconstruction of algebraic geometry. Those generalizations are the notion of scheme and topos. These two concepts are central in Grothendieck's work, and particularly relevant if one wants to understand Grothendieck's idea of space. In \cite{grot1} Grothendieck examined his more important contributions in mathematics and regarding toposes and schemes he writes: \\

``{\it Parmi ces th\`{e}mes, le plus vaste par sa portée me para\^{i}t \^{e}tre celui des topos, qui fournit l'idée d'une synth\`{e}se de la géométrie algébrique, de la topologie et de l'arithmétique. Le plus vaste par l'étendue des developpements auxquels il a donné lieu d\`{e}s \`{a} present, est le th\`{e}me des schémas... C'est lui qui fournit le cadre ``par excellence'' de huit autres parmi ces
th\`{e}mes envisagés, en m\^{e}me temps qu'il fournit la notion centrale pour un renouvellement de fond en comble de la géométrie algébrique, et du langage algébrico-géométrique.}'' \\

It is clear from this quote the unifying role played by schemes and toposes inside Grothendieck's mathematical vision. Thus, these two concepts allow to study not only topological or geometric aspects of the space but also algebraic and arithmetic ones. Here, we would like to emphasize that from the very beginning Grothendieck's idea of space involve a particular entity of arithmetic and geometric nature. That is what Zalamea \cite{zal} calls the space-number\footnote{In \cite{cruz1} we associate Zalamea's idea of the space-number to the idea of form in Grothendieck. We will go back to this point later in this text, since it is important in this discussion about the notion of space.} in analogy with the space-time in physics. We agree with Zalamea's interpretation, so we believe we need to understand Grothendieck's idea of space in terms of his idea of merging arithmetic and topology/geometry\footnote{This is a manifestation of the duality discrete/continuum. When we discuss Grothendieck's ideas of the space coming from his work in the 80's, this dichotomy will be clearer. Now, it is interesting to notice that some mixed discrete/continuum structure for the space is present in Grothendieck's first works in algebraic geometry, at least in an implicit way, and it will be more relevant when Grothendieck tries to think in the physical space.}. Therefore, schemes and toposes can be seen as the first instances of the space-number that Grothendieck had in mind.  

\subsection{Schemes: generalizing algebraic varieties} 

In his ICM talk in 1958 \cite{grot3} Grothendieck introduced his vast program of generalization of algebraic geometry that took him the next decade. During his years at the IHES he dedicated his life\footnote{Grothendieck's capacity of work is legendary. It is said that he could spend until 12 hours of continuous work.} to develop these ideas. One of the central notion he introduced is the notion of scheme and for us this is important since it is his first generalization of the notion of space. \\

As it is discussed in \cite{del1, zal} Grothendieck's intrepid idea is to accept that for every commutative ring $A$ (with unit) it is possible to define an affine scheme $X = \mathrm{Spec}(A)$. What is needed here is to define the structure sheaf $\tilde{A}$. In this case, the structure sheaf is the sheaf of rings over $X$ with a basis of open sets $\{ x \in X: f(x)= [f]_x \neq 0 \}$ with $f \in A$ and with fibres the local rings $A_x$ with $x \in X$. Thus, considering a ringed space, i.e. a pair $(X, \mathcal{O}_X)$ with $X$ a topological space and $\mathcal{O}_X$ a sheaf of rings over $X$, we say that it is an affine scheme if it is isomorphic to $(\mathrm{Spec}(A),\tilde{A})$. A general scheme is defined by gluing affine schemes. \\

Note that in the definition above the points of a scheme are not relevant. The important thing is that a scheme is an object where we can make {\it localization}, in the sense we discussed before. Thus, a scheme is a space of algebraic geometric nature that generalizes the notion of algebraic variety. Indeed, Grothendieck characterized the notion of scheme as a metamorphosis of the notion of algebraic variety (see \cite{grot1}). The use of the word metamorphosis is interesting, since it suggests a natural transformation\footnote{Naturalness is important for Grothendieck. In his view, mathematical concepts must be natural, in contrast to artificial and intricate constructions.} from the notion of algebraic variety to that of scheme.  \\

About naturalness of the notion of scheme, Grothendieck wrote in  \cite{grot1} \\

``{\it La notion de schéma est la plus naturelle, la plus``évident'' imaginable, pour englober
en une notion unique la série infinie de notions de ``variété'' (algébrique) qu'on maniait précédemment ...}'' \\

And also in the same text adds \\

``{\it L'idée m\^{e}me de schéma est d'une simplicité enfantine - si simple, si humble, que
personne avant moi n'avait songe \`{a} se pencher si bas. Si “béb\^{e}te” m\^{e}me, pour tout dire,
que pendant des années encore et en dépit de l'évidence, pour beaucoup de mes savants coll\`{e}gues, \c{c}a faisait vraiment ``pas sérieux''! Il m'a fallu d'ailleurs des mois de travail serré et solitaire, pour me convaincre dans mon coin que ``\c{c}a marchai'' bel et bien...}'' \\

It is clear that the notion of scheme is for Grothendieck the natural way in which the notion of algebraic variety should evolve, but also he noticed that the notion was always there, but it was needed to have the right eyes (his eyes in this case) to see it. Grothendieck was interested in capturing in one notion the multitude of varieties that arise when we work over prime characteristics. Again, the dichotomy discrete/continuum is present here. This is the fundamental aporia associated to the notion of space and Grothendieck started to struggling with it by introducing the idea of scheme. However, the notion of space is richer than that of algebraic variety and Grothendieck knew that.

\subsection{Toposes: spaces without points}

At this point we have just discussed a restricted idea of space since we have considered only those spaces that are algebraic varieties and look at their metamorphosis in the idea of scheme. However, for Grothendieck the space is something more general and, in this sense, the notion of scheme is not enough to capture the subtleties of space. Thus, it is necessary to go beyond algebraic geometry and entering in the realm of topology. The right notion to start doing that is the notion of topos. \\

In order to understand what a topos is we need to start with the idea of site. Here we will follow the description in \cite{zal}. Let us consider a category $\mathcal{C}$ and let us define a topology on it. For any object $X$ of $\mathcal{C}$ it is possible to define a set of sieves $J(X)$\footnote{A sieve can be identified with subobjects of the dual category.} satisfying some axioms, namely : \\

1. For any object $X \in \mathcal{C}$, any sieve $R \in J(X)$ and any morphism $f: Y \rightarrow X \in \mathcal{C}$, the sieve $R \times_X Y $ is in $J(Y)$. \\

2. For any sieves of $X$, $R$ and $R'$, if $R \in J(X)$ and for any $Y \rightarrow R$, $R' \times_X Y \in J(Y)$, then $R' \in J(X)$. \\

3. For any object $X$ of $\mathbb{C}$, $X \in J(X)$.  \\

A category $\mathcal{C}$ endowed with a topology that satisfies 1,2 and 3 is called a site. We focus on the idea of site since it allows us to illustrate Deligne's point that a space is something where localization makes sense. However, a site is just the first germ in the generalization of topological spaces that Grothendieck had in mind. Then, what is needed to be considered is a sheaf over a site. \\

Thus, a topos\footnote{In this text, we use the word topos in order to refer Grothendieck topos. There is a more general definition of topos (elementary topos) due to Lawvere but we do not consider it here. } is the category of all sheaves on a site. Therefore, a topos can be seen as a generalized space. Like schemes that can be seen as a metamorphosis of the idea of algebraic variety, a topos can be seen as a metamorphosis of the idea of topological space. Again, the word metamorphosis is quite accurate, since it allows to see how toposes arise as a natural deformation of the classical idea of topological space\footnote{We will see later that this metamorphosis does not stop in the notion of topos and will take another forms in terms of moderate topology and geometry of the forms.}. \\

Grothendieck introduced the idea in order to deal with cohomology theories in algebraic geometry, and in particular, toposes play an important role in Grothendieck's  approach to Weil conjectures. This way, toposes also arise as an important objects in problems of arithmetic nature. Again, it is possible to see that for Grothendieck the notion of space is something beyond purely geometric/topological considerations and it is central in whole mathematics. In \cite{grot1}, Grothendieck points out: \\

``{\it Celle de topos constitue une extension insoup\c{c}onnée, pour mieux dire, une métamorphose de la notion d'espace. Par l\`{a}, elle porte la promesse d'un renouvellement semblable de la topologie, et au del\`{a} de celle-ci, de la géométrie. D\`{e}s \`{a} présent d'ailleurs, elle a joué un r\^{o}le crucial dans l'essor de la géométrie nouvelle (surtout \`{a} travers les th\`{e}mes cohomologiques $\ell$-adique et cristallin qui en sont issus, et \`{a} travers eux, dans la démonstration des conjectures de Weil)}''.  \\

In the same text Grothendieck adds: \\

``{\it C'est le th\`{e}me du topos, et non celui des schémas, qui est ce ``lit'', ou cette ``rivi\`{e}re profonde'', o\`{u} viennent s'épouser la géométrie et l'alg\`{e}bre, la topologie et l'arithmétique, la logique mathématique et la théorie des catégories, le monde du continu et celui des structures ``discontinues'' ou ``discr\`{e}tes''. Si le th\`{e}me des schémas est comme le coeur de la géométrie nouvelle, le th\`{e}me du topos en est l'enveloppe, ou la demeure. Il est ce que j'ai con``\c{c}u de plus vaste, pour saisir avec finesse, par un m\^{e}me langage riche an résonances géométriques, une ``essence'' commune \`{a} des situations des plus éloignées les unes des autres, provenant de telle région ou de telle autre du vaste univers des choses mathématiques}.'' \\

From the last quote is very clear that for Grothendieck the notion of space, with its avatar in the notion of topos, is ubiquitous and central in mathematics \footnote{The notion of topos is, in certain sense, marginal in mathematics nowadays. However, recent works of Alain Connes and his collaborators (see \cite{conn}, as an example) have brought new life to the idea of topos and we believe this concept will play a central role in major mathematical problems like the Riemann hypothesis, for instance.} and adequate device to study in a unified way arithmetic and geometry. A place where discrete and continuous structures can live together in harmony.  \\

So far, we have seen a generalization of the notion of space running in two steps. First one corresponding to the notion of scheme and the second one corresponding to the notion of topos. This two steps generalization occurs inside the domain of algebraic geometry which was the main object of Grothendieck's reflections between 1958 and 1970. Thus, starting with the natural spaces in algebraic geometry, i.e. algebraic varieties, Grothendieck tried to find an archetypical structure for the space in the notion of topos. However, as we will see below the notion of topos does not constitute the final step in Grothendieck's thoughts related to the notion of space. There are other aspects to be considered and the idea of going beyond the notion of topos will dominate Grothendieck's work during the 80's. 

\section{Searching for new foundations: tame topology}

Around the middle of the 80's Grothendieck wrote an interesting text \cite{grot2} detailing several research programs he was thinking during the 70's, after his retirement of the mathematical milieu, that he considered would play an important role in the future mathematical research. One of these programs concerns to the development of the so-called tame topology, a new foundations for topology. In Grothendieck's words (see \cite{grot2}): \\

``{\it I would like to say a few words now about some topological considerations which have made me understand the necessity of new foundations for ``geometric'' topology, in a direction quite different from the notion of topos, and actually independent of the needs of so-called ``abstract'' algebraic geometry (over general base field and rings). The problem I started from, which already began to intrigue me some fifteen years ago, was that of defining a theory of ``dévissage'' for stratified structures, in order to rebuild them, via a canonical process, out of ``building blocks'' canonically deduced from the original structure}.'' \\

Here, in this quote, we see how Grothendieck wanted to propose a new setting that goes beyond the idea of topos and, in fact, beyond algebraic geometry. As we have seen, the origin and the development of the notion of topos is linked to the development of algebraic geometry in the setting Grothendieck envisioned during the 60's, so the notion of space developed during this period. Nevertheless, Grothendieck realized that in order to capture the essence of the space is necessary to considered a more general picture. Thus, with his idea of tame topology, the notion of space in Grothendieck started to be independent of the algebraic-geometric origin of Grothendieck's first approach. \\

Grothendieck chose an axiomatic approach in order to understand what a tame space could be. In his words: \\

``{\it My approach toward possible foundations for a tame topology has been an axiomatic one. Rather than  declaring (which would indeed be a perfectly sensible thing to do) that the desired ``tame spaces'' are no other than (say) Hironaka's semianalytic spaces, and then developing in this context the
toolbox of constructions and notions which are familiar from topology, supplemented with those which had not been developed up to now, for that very reason, I preferred to work on extracting which exactly, among the geometrical properties of the semianalytic sets in a space $\mathbb{R}^n$, make it possible to
use these as local ``models'' for a notion of ``tame space'' (here semianalytic), and what (hopefully!) makes this notion flexible enough to use it effectively as the fundamental notion for a ``tame topology'' which would express with ease the topological intuition of shapes}. '' \\

Even though, as Grothendieck wrote, a tame space is not a semianalytic space, it is illustrative, in order to understand the idea Grothendieck had in mind, to discuss a bit the notion of semianalytic set, since Grothendieck did not provide an explicit definition of tame space in \cite{grot2}. \\

Let $M$ be a real analytic manifold and consider a subset $X$ of $M$, the $X$ is called semianalytic if each $a \in M$ has a neighbourhood $\mathcal{U}$ such that $X \cap \mathcal{U} \in S(\mathcal{O}(\mathcal{U}))$, where $\mathcal{O}(\mathcal{U})$ denotes the ring of real analytic functions on $\mathcal{U}$ and $S(\mathcal{O}(\mathcal{U}))$ denotes the smallest family of subsets of $\mathcal{U}$ containing all $\{ f(x) >0 \}$, $f \in  \mathcal{O}(\mathcal{U})$, which is stable under finite intersections, finite unions and complement. \\

Without a precise definition, but using the notion of semianalytic space as local model, it is possible to see that the new tame spaces fit in the general definition of space we have assumed following Deligne. However, there is an important point in the origin of the idea of tame spaces that it did not appear in the previous notions of scheme and topos. When Grothendieck formulated the idea of tame space he was thinking in a right framework to deal with the notion of geometrical shape. Thus, the main motivation here is to develop a notion of space that allows us to study the geometric properties of shapes\footnote{In the next section we will see that with the same motivation of understanding the elusive idea of form (being geometrical shapes just one instance of this archetypical notion)  Grothendieck was beyond tame spaces and proposed the so-called geometry of forms.}. This point is very clear from this quote (see \cite{grot3}): \\

`` {\it After some ten years, I would now say, with hindsight, that ``general topology'' was developed (during the thirties and forties) by analysts and in order to meet the needs of analysis, not for topology per se, i.e. the study of the topological properties of the various geometrical shapes. That the foundations of topology are inadequate is manifest from the very beginning, in the form of ``false problems'' (at least from the point of view of the topological intuition of shapes) such as the ``invariance of domains'', even if the solution to this problem by Brouwer led him to introduce new geometrical ideas}.'' \\

Unlike the notion of scheme and topos, both central in algebraic geometry and widely known for mathematicians nowadays, the notion of tame space, and therefore tame topology, is marginal in mathematics. There are some developments in logic \cite{van} in relation to O-minimal structures, but besides this the idea of tame space is not that relevant\footnote{We consider that it would be a very interesting research program to try to develop Grothendieck's ideas in \cite{grot2} related to tame topology. We think that this could open new geometric perspectives and new understanding of topological phenomena as Grothendieck envisioned.}. However, in our study of the notion of space in Grothendieck is, in fact, very important, since it provides the first attempt of Grothendieck to formulate a theory of space with the clear motivation of providing a setting to study shapes as foundational objects for geometry/topology, as we already mentioned. In this sense, we believe that it is not possible to understand how the idea of space evolved in Grothendieck's work without passing through the discussion about tame topology. On the other hand, we will see that in order to tackle the problem of understanding shapes, this is not the final answer that Grothendieck wanted to provide. \\

It is important to remark that the notion of tame space does not pretend to generalise the idea of topological space (like toposes do) but provides a new understanding for it. In that sense, a tame space is a new form that the notion of space can take, it is new type of space. Thus, we can see that the notion of space is really protean for Grothendieck. 

\section{Towards a new geometry of forms}

We saw in the last section how part of the motivation for developing a tame topology was to construct a right notion of space in order to deal with the various geometrical shapes. In fact, the notion of form\footnote{In \cite{cruz1} we have argued that there is an important difference between form and shape and we have discussed how the notion of form is more general than that of shape. In fact, Grothendieck himself chose to use the word form rather than shape, as we can see in a letter to Ronald Brown on June 17th 1986: \\
{\it ``I've been very strongly involved with mathematics lately, working out another approach to ``topology'' and ``form'' (different from topological spaces, from topoi and from moderate spaces as proposed in Esquisse d'un Programme), just getting the basic language straight ...'' }} was one of the main subjects of Grothendieck's reflections in the 80's, and in some sense tame topology was not enough for Grothendieck in order to provide a general setting for the notion of space if we want to study the geometric properties of forms. In this direction, Grothendieck proposed a new approach for the foundations of topology. This is what he called geometry of forms. In a letter to Yamashita on July 9th 1986, Grothendieck wrote: \\

``{\it I have been very intensely busy for about a month now, with writing down some altogether different foundations of ``topology'', starting with the ``geometrical objects'' or ``figures'', rather than with a set of ``points'' and some kind of notion of ``limit'' or (equivalently) ``neighbourhoods''. Like the language of topoi (and unlike the so-called ``moderate space'' theory foreshadowed in the Esquisse, still waiting for someone to take hold of the work in store ...), it is a kind of topology ``without points'' - a direct approach to ``shape''. I do hope the language I have started developing will be appropriate for dealing with finite spaces, which come off very poorly in ``general topology'' (even when working with non-Hausdorff spaces)}'' \\

Grothendieck wrote a manuscript \cite{mont} about these ideas. The manuscript remained relatively unknown (except for a few mathematicians) until the publication online of Grothendieck's manuscripts by Université de Montpellier. However, due to Grothendieck's handwritten it is not easy to decipher its content. To the best of our knowledge the first places where there is an attempt to discuss the content of the manuscript are \cite{cruz2, cruz3} based on some correspondence of Grothendieck where he wrote about it. For the purpose of this essay, we consider that the ideas on the manuscript that we can extract from the correspondence are very important since they allow us to really understand how the notion of space was evolving in Grothendieck's work. We want to argue that a full understanding of the problem of space in Grothendieck is impossible without studying the manuscript or at least its main ideas. However, a deep study of the text is still a work in progress, but here we hope to give a good idea of what Grothendieck was thinking\footnote{Just to provide an idea of the text, we want to give the title of the different chapters of the manuscript. \\
Chap 1. On a topology of the (topological) forms.
Chap 2. Topological realisations of networks. 
Chap 3. Networks via decompositions. 
Chap 4. Analysis situs (First attempt).
Chap 5. Algebra of figures. 
Chap 6. Analysis situs (Second attempt). 
Chap 7. Analysis situs (Third attempt). 
Chap 8. Analysis situs (Fourth attempt). 
Chap 9. Notes.  \\
It is interesting to see how Grothendiek uses the expression analysis situs that evokes Leibniz and Poincaré. In a letter to Yamashita on September 16th 1986, Grothendieck mentioned Poincaré conjecture in relation to the analysis situs he was invented. It is still matter of study to understand the relation between Poincaré and Grothendieck (and also between Grothendieck's ideas and Leibniz' ideas), but we find illustrative this extract of the letter: \\ 
``{\it I have been impressed, of course, by your second list of ``newest events''. I didn't even know about the solution of Poincar\'e's conjecture - comes just at the right time for me, to be able to justify within the framework of the ``geometry of forms'' or ``analysis situs'' I am developing in terms of extant ``general topology'', a certain definition of ``regular figure'' (the combinatorial substitute for ``variety'') I had in mind...''}. \\
Grothendieck refers the proof announced by Colin Rourke and Eduardo R\^{e}go.  The proof turned out to be incorrect. However, some years after this, Gregory Perelman found a right proof for the conjecture. So it is tempting to ask: What is the relation between Poincaré conjecture and Grothendieck's geometry of forms?}. \\

In the extract of the letter to Yamashita we have quoted there is one important point that we want to remark. Grothendieck pointed out that his new approach to the notion of space is closer to the idea of topos than the idea of tame space, in the sense that these new spaces are not made by points, but also because here he is looking for a direct approach to shapes. However, looking at the discussion in \cite{grot3} concerning tame topology we can see that the motivation of Grothendieck was to understand the forms from a geometric point of view. Is there any contradiction here? We think the answer is no. \\

One thing is to be motivated by some idea and another thing is that this idea can be carried out by certain construction we propose. Our conjecture is that Grothendieck realized that his tame topology was not the right framework for the goal of understanding forms, taking forms in a very broad sense, and for this reason he saw the necessity  of producing a new theory and that is what he expressed in the letter to Yamashita. Actually, Grothendieck explored different approaches during the years in order to understand what exactly a space is and the idea of tame space can be thought as one of them. Before his retirement at Lasserre in June 24th 1991, Grothendieck wrote a letter to an unknown recipient\footnote{The unknown recipient is identified as A.Y. This letter and its relation with the manuscript about the geometry of forms has been discussed in \cite{cruz2}.}:\\

``{\it Il est vrai qu'us cours de dix derni\`eres ann\'ees, j'ai r\'efl\'echi ici et l\`a \`a diverses extensions de la notion d'espace, en gardant \`a l'espirit la remarque p\'en\'etrante de Riemann. J'en parle dans quelques lettres \`a des amis physiciens ou ``relativistes''. Il ne doit pas \^etre tr\`es difficile p. ex. de d\'eveloper une sorte de calcul diff\'erentiel sur des ``vari\'et\'es'' qui seraient des ensembles finis (mais \`a cardinal ``tr\`es grand''), ou plus g\'en\'eralement discretes, visualis\'es comme format une sorte de ``r\'eseau'' tr\'es serr\'e de points dans una vari\'et\'e $C^{\infty}$ (p. ex. una vari\'et\'e riemannienne) - une sorte de g\'eometri\'e diff\'erentielle ``floue'', o\'u toutes les notions num\'eriques son d\'efines seulement ``\`a $\epsilon$ pr\'es'', pour un ordre d'approximation $\epsilon$ donn\'e}.'' \\ 

This quote is particularly relevant. Three important points appear here that are not present in the reflections about the notion of space that Grothendieck proposed in his work in algebraic geometry and in tame topology. First of all, Grothendieck asks for the intrinsic structure of the space and its discrete or continuous structure. Of course, when he discussed about toposes the relation between the discrete and the continuum world was present, but in the sense on how toposes can put arithmetic and geometric phenomena together, but here the difference is that the discussion on the dichotomy continuum/discrete appears in a more foundational and structural way. Secondly, he was not only thinking of the mathematical space but also in the physical space. It is widely believed that Grothendieck was not interested in physics at all, but it seems that in his last period and regarding the problem of space he was interested in physics and, in fact, he had contact with some physicists and discussed with them his ideas. Finally, he recognise explicitly that some of his thoughts about the space came from Riemann. We will back to the relation with Riemann later, so let us focus on the first and the second points we have mentioned. \\

What was the kind of space that Grothendieck had in mind? We can find a first answer in the letter to Yamashita on June 9th 1986. There Grothendieck argued: \\

``{\it As Riemann pointed out I believe, the mathematical continuum is a convenient fiction for dealing with physical phenomena, and the mathematics of infinity are just a way of approximating (by simplification through ``idealisation'') an understanding of finite aggregates, whose structure seems too elusive or too hopelessly intricate for a more direct understanding (at least it has been so till now). Yet it may well be that we are approaching at present the point where the continuous models of the physical world fail to be adequate - but the physicists are so accustomed taking for granted the conceptual superstructure of ``continuous space'' worked out by the Greeks and their successors up to Relativity and Quantum mechanics, and to confuse it with reality, that they may well never get aware of a widening gap between concepts and phenomena. This brings to my mind that not only the language of (so-called ``topological'') shape is to be rethought from scratch, in order to be capable to account for ``shape'' of finite spaces, but the same holds with vast sectors of geometry and analysis, if not all of it. Thus differential and integral calculus, differential equations and the like, differentiable or Riemannian varieties, the tensor calculus on it, etc - all this quite evidently (once you start thinking about it) should make and has to make sense within the framework of (say) finite aggregates as well. Just think, say, of $10^{100}$ or $10^{1000}$ points pretty densely located on the unit sphere or in the unit ball, a lot more than needed for exceeding by very far he accuracy of any conceivable physical measurement, so that no physicists nor anybody whosoever (except maybe God Himself) could ever possible distinguish between the ambient ``continuous'' space (granting it does have some kind of ``physical'' existence, which I doubt...) and this kind of network thought of as an ``approximation'' (whereas in reality, the opposite seems to me more likely to be true: the continuum is the concept approximating the elusive finite - but very large - aggregate). This would mean rethinking too the notion of differentiability (say) and of differential, in wholly new terms, when ``passing to the limit'' means now a finite process, namely passing to ``small'' values of the parameter or parameters, where the rate of ``smallness'' is precisely prescribed and should by all means not be exceeded, and where the ``margin of errors'' admitted in the computation or definition of the ``differential'' or ``derivative'' is equally precisely prescribed, in keeping with the former prescription}.'' \\ 

Grothendieck was inclined for a discrete structure for the space and even he thought of the continuum as a convenient fiction. This is actually quite impressive. Nobody before Grothendieck proposed to think the continuum as an approximation to the discrete\footnote{In the letter Grothendieck attributes this idea to Riemann. However, it is not that clear that actually Riemann was thinking the continuum in this terms. We will back to this point later.}. We find this idea very suggestive on its own\footnote{We discuss about this in \cite{cruz2}. Let us just mention that this idea is an inversion respect to the usual way of thinking the relation between the continuum and the discrete. Usually the continuum is seen approximate by the discrete via a saturation (limit) procedure. Grothendieck is proposing a radical, but we think very fruitful, change of perspective. It would be very interesting to explore in more detail this line of thought.}. For Grothendieck, the continuous structure of the space can be thought as an artificial, but necessary and convenient, construction in order to understand its actual structure. However, the relation between the discrete and the continuum is not that simple for Grothendieck. In fact, in \cite{grot1} he proposed that there could be three possibilities for the structure of the space. The space can be discrete, continuum or mixed discrete/continuum structure. In Grothendieck's words: \\

``{\it Les développements en mathématique des derni\`{e}res décennies ont d'ailleurs montré une symbiose bien plus intime entre structures continues et discontinues, qu'un ne l'imaginait encore dans la premi\`{e}re moitié de ce si\`{e}cle. Toujours est-il que de trouver un mod\`{e}le ``satisfaisant'' (ou, au besoin, un ensemble de tels mod\`{e}les, se ``raccordant'' de fa``\c{c}on aussi satisfaisante que possible. . .), que celui-ci soit ``continu'', ``discret'' ou de nature ``mixte'' — un tel travail mettra en jeu s\^{u}rement une grande imagination conceptuelle, et un flair consommé pour appréhender et mettre \`{a} jour des structures mathématiques de type nouveau. Ce genre d'imagination ou de ``flair'' me semble chose rare, non seulement parmi les physiciens (o\`{u} Einstein et Schr\"{o}dinger semblent avoirété parmi les rares exceptions), mais m\^{e}me parmi les mathématiciens (et l\`{a} je parle en pleine connaissance de cause)}.'' \\

Again, it is clear, that he was thinking not just in the mathematical space but also in the physical one. Unlike what he did in algebraic geometry and in tame topology, his new geometry of forms is not just a mathematical description of the space but it should be thought as a description useful for physical purposes. This is the reason, we believe, that the notion of form was so central in this new approach. A form has many manifestations in the world not only in mathematics, so any theory about the space focused on the idea of form should reflect those manifestations. We think Grothendieck had fully conscience about this. 

\subsection{Riemann's heritage?} 

An important point to analyse in Grothendieck's conception of the geometry of forms is Riemann's influence. Of course, Riemann has a big influence in the whole grothendieckean work, but in relation to the geometry of forms that influence is particularly interesting. In \cite{grot1} we find: \\

``{\it Il doit y avoir déj\`{a} quinze ou vingt ans, en feuilletant le modeste volume constituant l'oeuvre compl\`{e}te de Riemann, j'avais été frappé par une remarque de lui ``en passant''. Il y fait observer qu'il se pourrait bien que la structure ultime de l'espace soit ``discr\`{e}te'', et que les représentations ``continues'' que nous nous en faisons constituent peut-\^{e}tre une simplification (excessive peut-\^{e}tre, \`{a} la longue. . .) d'une réalité plus complexe; que pour l'esprit humain, ``le continu'' était plus aisé \`{a} saisir que ``le discontinu'', et qu'il nous sert, par suite, comme un ``approximation'' pour appréhender le discontinu. C'est l\`{a} une remarque d'une pénétration surprenante dans la bouche d'un mathématicien, \`{a} un moment o\`{u} le mod\`{e}le euclidien de l'espace physique n'avait jamais encore été mis en cause; au sens strictement logique, c'est plut\^{o}t le discontinu qui, traditionnellement, a servi comme mode d'approche technique vers le continu}.'' \\

Here Grothendieck again refers Riemann as the one who was the first in thinking the continuum as an approximation of the discrete, and in this sense the continuous structure of the space as an approximation of a more intricate discrete structure of it\footnote{We already saw that Grothendieck also mentions Riemann in relation with this idea in his letter to Yamashita on June 9th 1986}. What it is also important to consider is that Grothendieck pointed out that the relation between the continuum and the discrete can be determined by the fact that the continuum is easier to assimilate by the human mind (he uses the expression human spirit) that the discrete. This has deep philosophical implications in terms of the nature of space. The question is to some extent this idea is already in Riemann in the explicit way that Grothendieck is presenting it. \\

Riemmann's contribution started to be clarified in a letter to Yamashita on September 16th 1986. There Grotehndieck wrote:  \\

``{\it When I wrote about the theory of mathematical continuum begin a kind of approximation and idealization of very large and complex discrete aggregates, I thought I was citing an idea of Riemann, from what I remembered from reading through his collected works, many years ago (presumably in his exposition of the idea of what is now called riemannian geometry). When I later looked for a precise reference in his work, I was rather surprised that I wasn't able to find it - and it seems you didn't find it either. So maybe, I found his idea only ``in between the lines'' of his actual words}. '' \\ 

The ``in between the lines'' at the end is really remarkable. At the end of his Habilitation thesis \cite{riem}, Riemann discussed about the possibility of a discrete structure for the space. That idea yielded a big impression on Grothendieck. After looking with more detail, he recognised that in an explicit way the idea is not in Riemann but he certainly believed that it can be extracted from what Riemann wrote. This is a risky implication. Thus, Why does Grothendieck make that implication? We could find a partial answer again in the letter to Yamashita from September 16th: \\

``{\it It would almost seem that there has been something like ``telepathic'' communication (to give it a name) from Riemann's thought-world to mine. I am still convinced now, or rather, I well know, that I never discovered myself this idea (because this kind of discovery of some significant view, I mean the moment of such discovery, one never forgets I believe) but that I got it from Riemann, and with some surprise...}'' \\

This passage is very beautiful. That ``telepathic communication'' between Grothendieck and Riemann, two of the biggest mathematical geniuses of all times, connected by their mathematical passion and the wish of understanding the notion of space, is one of the most remarkable examples of mathematical creativity.  Despite of being separated by many years, both of them had the intention of unveiling the inner and deep structure of the space and they walked together going behind that goal. Thus, without doubt Grothedieck's ideas about the space, and in particular the geometry of forms he envisioned in the late 80's, is Riemann's heritage.  

\section{The notion of space: mathematics, physics and philosophy}

In the previous section we discussed how thanks to the formulation of his geometry of forms Grothendieck began interested in the relations between the physical and the mathematical space. It is also important to remark that in fact Grothedieck was beyond mathematics and physics in his reflections and also consider the importance of philosophy in order to understand what a space is. This give us a new perspective about Grothedieck's thoughts in the late 80's and also about his understanding of the notion of space as a whole. In particular, Grothendieck was interested in the role of mathematical models for the physical reality and the philosophical background of them. He thought that only putting together philosophy, physics and mathematics we can find an adequate approach to the problem of space. Regarding this, in \cite{grot1} he wrote: \\

``{\it Pour résumer, je prévois que le renouvellement attendu (s'il doit encore venir. . .) viendra, plut\^{o}t d'un mathématicien dans l'\^{a}me, bien informé des grands probl\`{e}mes de la physique, que d'un physicien, Mais surtout, il y faudra un homme ayant ``l'ouverture philosophique'' pour saisir le noeud du probl\`{e}me. Celui-ci n'est nullement de nature technique, mais bien un probl\`{e}me fondamental de ``philosophie de la nature''}.'' \\

And in a letter to AY on June 24th 1991 he also pointed out: \\

``{\it ``Je suspecte que les nouvelles structures \`a d\'egager seront beaucoup plus subtiles qu'un simple paraphrase de mod\`eles continus connus en termes discretes. Et surtout, qu`avant toute tentaive de d\'egager des nouveaux mod\`eles, pr\'esm\'es meilleurs que les anciens, il s\'impose de poursuivre une r\'eflexion philosophico-math\'ematique tr\`es servie sur la notion m\^eme de ``mod\`ele math\'ematique'' de quelque aspect de la realit\'e - sur son r\^ole, son utilit\'e, et sus limites}.'' \\

It is beautiful to see how Grothendieck uses the expression natural philosophy. The problem of space is not just a mathematical, physical of philosophical one, it is a problem of natural philosophy in a very broad sense that includes mathematics, physics and philosophy but also goes beyond them. This is evidence of how deep Grothendieck's thoughts about the space became. The philosophical side of Grothendieck's approach to the notion of space is, in general, unknown but for us it is very relevant and it should be studied in more detail. Here, we just mention this aspect to motivate the study and the debate about it. \\

Grothendieck pointed out that the problem of space is not only of technical nature. It is not enough to think of the problem in mathematical or physical terms. We also need a right conceptual approach and it is in this point where philosophy should enter. One of the aspects to be consider is how the mathematical models of the space that we can propose reflect the nature and the structure of the physical space and it is consistent philosophically. Of course, this is not an easy problem and Grothendieck had no answer. Again, the reflection on the role of mathematical models is of philosophical nature as Grothendieck saw (see \cite{grot1}): \\

``{\it  J’ai le sentiment que la r´eflexion fondamentale qui attend d'\^{e}tre entreprise, aura \`{a} se placer sur deux niveaux différents. \\

1. Une réflexion de nature ``philosophique'', sur la notion m\^{e}me de ``mod\`{e}le mathématique'' pour une portion de la réalité. Depuis les succ\`{e}s de la théorie newtonienne, c'est devenu un axiome tacite du physicien qu'il existe un mod\`{e}le mathématique (voire m\^{e}me, un mod\`{e}le unique, ou ``le'' mod\`{e}le) pour exprimer la réalité physique de fa\c{c}on parfaite, sans ``d'ecollement'' ni bavure. Ce consensus, qui fait loi depuis plus de deux si\`{e}cles, est comme une sorte de vestige fossile de la vivante vision d'un Pythagore que ``Tout est nombre''. Peut-\^{e}tre est-ce l\`{a} le nouveau ``cercle invisible'', qui a remplacé les anciens cercles métaphysiques pour limiter l'Univers du physicien (alors que la race des ``philosophes de la nature'' semble définitivement éteinte, supplantée haut-la-main par celle des ordinateurs. . .). Pour peu qu'un veuille bien s'y arr\^{e}ter ne fut-ce qu'un instant, il est bien clair pourtant que la validité de ce consensus-l\`{a} n'a rien d'évident. Il y a m\^{e}me des raisons philosophiques tr\`{e}s sérieuses, qui conduisent `a le mettre en doute a priori, ou du moins,\`{a} prévoir \`{a} sa validité des limites tr\`{e}s strictes. Ce serait le moment ou jamais de soumettre cet axiome \`{a} une critique serrée, et peut-\^{e}tre m\^{e}me, de ``d'emontrer'', au del\`{a} de tout doute possible, qu'il n'est pas fondé : qu'il n'existe pas de mod\`{e}le mathématique rigoureux unique, rendant compte de l'ensemble des phénom\`{e}nes dits ``physiques'' répertoriés jusqu'\`{a} présent. \\

Une fois cernée de fa\c{c}on satisfaisante la notion m\^{e}me de “mod\`{e}le mathematique”, et celle de la “validite” d'un tel mod\`{e}le (dans la limite de telles ``marges d'erreur'' admises dans les mesures faites), la question d'une ``théorie unitaire'' ou tout au moins celle d'un ``mod\`{e}le optimum'' (en un sens \`{a} préciser) se trouvera enfin clairement posée. En m\^{e}me temps, on aura sans doute une idée plus claire aussi du degré d'arbitraire qui est attaché (par nécessité, peut-\^{e}tre) au choix d'un tel mod\`{e}le. \\

2. C'est apr`es une telle réflexion seulement, il me semble, que la question ``technique'' de dégager un
mod\`{e}le explicite, plus satisfaisant que ses devanciers, prend tout son sens. Ce serait le moment alors, peut-\^{e}tre, de se dégager d'un deuxi\`{e}me axiome tacite du physicien, remontant \`{a} l'antiquité, lui, et profondément ancré dans notre mode de perception m\^{e}me de l'espace : c'est celui de la nature continue de l'espace et du temps (ou de l'espace-temps), du ``lieu'' donc o\`{u} se déroulent les ``phénom\`{e}nes physiques''}.''\\

We include this extensive quote because it illustrates in Grothendieck's own words what we are trying to argue, i.e., that the technical description of the space only makes sense for Grothendieck after a conceptual approach. It took him around 40 years to reach this conclusion. This quote also shows the evolution of Grothendieck's ideas from his pure technical questions in his work in algebraic geometry to the more physical-philosophical questions in his geometry of forms. Thanks to this, it is very clear that notion of space is indeed a leitmotif in Grothedieck's mathematical life. 

\section{Final remarks}

Along this essay we have shown how the idea of space in Grothedieck evolved from his generalization of the notion of (algebraic) variety via schemes to the formulation of a new approach for the foundations of topology in what he called geometry of forms, passing through the notion of topos and tame space. Of course, for Grothendieck the notion of variety and topological space (in the classical sense) were useful in mathematics but he realized of new situations where new ideas on space were needed. In particular, one of these situations is related to the notion of form and its role in topology and, in general, in mathematics. \\

Grothendieck wrote that there are three main areas of mathematical reflection, namely arithmetic (associated to the notion of integer numbers), analysis (associated to the notion of magnitude) and geometry (associated to the notion of form). From the beginning he put the notion of form as belonging to the realm of geometry and more general topology, but even more, it is the central notion for geometry in the sense that geometry for him is the branch of mathematics in charge of studying the aspects of the universe, not just the mathematical universe of course, related to the form.\\

Following that line, for Grothendieck there are three different notions which are important from a mathematical point of view, namely  number, magnitude and form. In \cite{cruz1} we argue that number and magnitude are manifestations of a more general notion of form and in this sense number can be thought as an arithmetical form while magnitude can be thought as an analytic form. What Grothendieck called form should be interpreted as a geometrical form (or shape). What we have is a fundamental triad that gives the most important aspects to be considered in a mathematical approach to the notion of space, i.e. arithmetical, analytical and geometrical aspects. We have argued that this was the quest of Grothendieck for more than 40 years. \\

On the other hand, as we have already mentioned before in this text, Zalamea's interpretation \cite{zal} of the space-number in Grothendieck's work as analogue to the space-time in physics gives a nice way to understand Grothendieck's ideas. Zalamea proposes that we should think Grothendieck's work as an attempt to develop the notion of space-number as the conjunction of arithmetic and geometry in one entity. This analogy is very suggestive since, as we have shown, Grothendieck was interested in physical aspects of the space. For this reason, we find very plausible to think of his last works in terms of constructing a theory for a space-number. We think that this line of thought could be very helpful in order to understand the notion of space in Grothendieck and it is a future project try to follow it. \\

In \cite{bol} Bolzano presents an idea of mathematics as a general theory of forms. In some sense Grothendieck's approach to the notion of space is a realization of Bolzano's idea, including physics and philosophy in the picture. We believe this should be a way to think of Grothendieck's contribution to the human knowledge. There are still many problems to be considered and many gaps to be filled in order to understand Grothendieck's ideas about the space. In this text we wanted to present a new perspective and contribute with some elements which, to the best of our knowledge, were not consider before in the literature. More than give a definitive answer about the problem of space in Grothendieck's work we wanted to give a picture, as complete as we could do it, in order to start the dialogue. A dialogue that following Grothendieck's own advice will require a mathematical, physical and philosophical way of thinking. \\ \\

{\bf Acknowledgements} The author wants to thank Fernando Zalamea for many illuminating discussions on Grothendieck's work over the years, his constant and warm support and having read a draft version of this essay. Of course, any mistake that remains is completely responsibility of the author. He also thanks Winfried Scharlau for providing him with Grothendieck's correspondence with Yamashita and Tsuji and for conversations on Grothendieck's life, and Yuri Manin for supporting the idea of writing this text. Finally, he thanks Pierre Deligne for many interesting emails and his patience in explaining his ideas.

\end{document}